\documentclass{article}
\usepackage[utf8]{inputenc}
\usepackage[margin=1in]{geometry}
\usepackage{amsmath} % assumes amsmathpackage installed
\usepackage{amssymb} % assumes amsmathpackage installed
\usepackage{mathtools}
\usepackage{acronym}
\usepackage{mathrsfs}

\usepackage{subcaption}
\usepackage{hyperref}
\hypersetup{colorlinks=true, linkcolor=black}
\usepackage{algorithm}%
\usepackage{algorithmicx}%
\usepackage{algpseudocode}%

\newcommand{\lure}{Lur\'{e} }

\newcommand{\R}{\mathbb{R}}

\newcommand{\N}{\mathbb{N}}

\newcommand{\Z}{\mathbb{Z}}

\newcommand{\0}{\mathbf{0}}
\newcommand{\1}{\mathbf{1}}

\DeclarePairedDelimiter{\abs}{\lvert}{\rvert}
\DeclarePairedDelimiter{\norm}{\lVert}{\rVert}

\DeclarePairedDelimiter{\diag}{\textrm{diag}(}{)}

\DeclarePairedDelimiterX{\inp}[2]{\langle}{\rangle}{#1, #2}

\DeclareMathOperator*{\find}{find}
\usepackage{amsthm}
\newtheorem{thm}{Theorem}
\newtheorem{prop}[thm]{Proposition}% 

\newtheorem{prob}{Problem}

\title{Incremental Peak-to-Peak Gain Computation and Regulation of Discrete-Time Positive Lur\'{e} Systems \\ using Linear Programming} % Title, preferably not more 
\author{Jared Miller    % Add the 
\renewcommand\footnotemark{}

\thanks{Chair of Mathematical Systems Theory, Department of Mathematics,  University of Stuttgart, Stuttgart, Germany. (email: \url{jared.miller@imng.uni-stuttgart.de})}         
\thanks{
This work was funded by Deutsche Forschungsgemeinschaft (DFG, German Research Foundation) under Germany's Excellence Strategy - EXC 2075 – 390740016. We acknowledge the support by the Stuttgart Center for Simulation Science (SimTech).
}
}

\date{\empty}
\begin{document}

\maketitle

\begin{abstract}
Incremental gains provide a quantitative measure of robustness for trajectories of a dynamical system.  Finite incremental gains ensure that all pairs of trajectories will together converge to a fixed point or will diverge together in the absence of an applied input.
This work approaches the problem of computing the incremental peak-to-peak gains for  positive linear systems in Lur\'e feedback with static memoryless nonlinearities, and regulating the incremental  $\ell_\infty$ gain through the design of a state-feedback controller. 
 Upper-bounds on these incremental gains can be computed through linear programming. Computation and state-feedback regulation of the  $\ell_\infty$ incremental gain is verified by numerical examples.

\end{abstract}
% \end{frontmatter}

\section{Introduction}
\label{sec:intro}

Positive systems {describe dynamics in which the inputs, outputs, and possibly the describing states obey structural nonnegativity properties \cite{farina2000positive}.}
These systems are of fundamental interest in mathematics, networks, biology, and planning \cite{brown2007compartmental}, \cite{shorten2006positive}.
Instances of positive systems  in control include management of inventory in supply chains, adjusting mixing rates in chemical batch processes,  {and fisheries management. The respective nonnegative quantity in each case is the number of items in a warehouse, concentration of a chemical, or number of fish in a lake.}
The robustness of these positive systems against uncertainties when performing tracking can be quantified by their incremental gains.
This work bounds the incremental $\ell_\infty$ (peak-to-peak) gain attained by the  \lure interconnection \cite{lur1944theory} of a positive discrete-time linear system and a static but possibly time-varying  nonlinearity.
Linear programming techniques can upper-bound the incremental $\ell_\infty$ gain of a given \lure-interconnected system using max-type lyapunov linear copositive functions. This method can be used to design state-feedback controllers that have a certificate of incremental $\ell_\infty$-bounded gains. These results are built upon the prior research in \cite{piengeon2024linear}  about incremental $\ell_1$ gain analysis under \lure bounded uncertainty.

{A system is called externally positive if every nonnegative input from zero initial state
results in a nonnegative output sequence, and is called internally positive if the state is
also required to be nonnegative for nonnegative inputs and {initial states}}. Internal positivity can be verified by checking signs of matrices in a state-space representation (nonnegative or Metzler). 
Testing external positivity of a state-space linear system is generically a co-NP hard problem \cite{blondel2002presence}, though convex heuristics exist to attempt recognition and construction \cite{grussler2017externalid}, \cite{grussler2021second}. The work in \cite{drummond2023externally}, \cite{taghavian2023external}, \cite{van2025externally} uses transfer function descriptions in order to recognize external positivity.

Stabilization of  internally positive systems can be accomplished using (dual) linear copositive Lyapunov functions, in which stabilizing controllers arise from the solution of efficiently solvable linear programming problems \cite{rantzer2015scalable} \cite{rantzer2018tutorial}. 
Switched positive systems can also be stabilized in this manner, by crafting {L}yapunov functions that are compatible between switching arcs \cite{fornasini2011stability} \cite{blanchini2015switched} \cite{FORNI20174558}. The performance criteria of the $\ell_\infty$ and $\ell_1$ norms can be modeled as Integral Linear Constraints,  and therefore can also be bounded  using linear programming \cite{ebihara20111} \cite{briat2013positive} \cite{ebihara2016analysis}. Linear programming method can also be used to synthesis state-feedback \cite{ebihara2012optimal}  and static-output feedback \cite{shen2016static} controllers with bounded performance.
We note that the work in \cite{rantzer2016kyppos} generates a KYP lemma for positive settings, which can be used to bound $\ell_2$ performance or verify passivity. The $H_2$ or $H_\infty$ norms of positive systems \cite{kato2019analysis} can be computed using NP-hard copositive optimization \cite{dur2010copositive}.

Incremental gains are closely related to contraction rates \cite{lohmiller1998contraction}. For general linear systems, the incremental and non-incremental gains are equivalent \cite{desoer2009feedback}. However, this equivalence of gains may break down when considering nonlinearities. 
In the positive systems setting, the increments of a positive system are not guaranteed to be positive, as positive systems form a conic behavior rather than a linear behavior \cite{padoan2023data}. {A bounded $\ell_\infty$-gain of a positive linear system under a nonlinear interconnection may  not imply that the incremental $\ell_\infty$ gain is also bounded.}

This paper has the following structure: Section \ref{sec:preliminaries} reviews preliminaries of notational { and  positive linear systems}. Section \ref{sec:gain_compute} presents linear programs to compute the incremental $\ell_1$ and $\ell_\infty$ gains. Section \ref{sec:regulate} extends the $\ell_\infty$ gain computation linear program for state-feedback stabilization. Section \ref{sec:examples} deploys this incremental bound on numerical examples. Section \ref{sec:conclusion} concludes the paper. Appendix \ref{app:proofs} contains all supporting proofs.

\section{Preliminaries}
\label{sec:preliminaries}
This section presents preliminaries of notation and positivity of linear systems.
\subsection{Notation}
\label{sec:notation}

The $n$-dimensional real vector space will be denoted as $\R^n$, and its nonnegative orthant and positive orthants will respectively be denoted as $\R^n_{\geq 0}$ and $\R^n_{>0}$. The set of $n \times m$ real-valued matrices is $\R^{n \times m}$. The transpose of a matrix $M \in \R^{n \times m}$ is $M^\top$. 

The notation $ \1, \0$ will respectively refer to {the} all ones and all zeros matrices. These matrices also include subscripted dimensions (e.g. $\1_n$ as a vector or  $\0_{m \times n}$ as a matrix) if not otherwise specified. The set of integers is $\Z$. The set of natural numbers {including 0} is $\N$. The subset of natural numbers between $a$ and $b$ inclusive is $a..b$. 
{The elementwise division of two vectors} $u, v \in \R^n$, {is} $u ./ v \  := \ \{u_i  / v_i\}_{i=1}^n$.

The elementwise absolute value of a vector (or matrix) $x$ is $\abs{x}$. The $L_p$ norm of a vector  is $\norm{x}_p$. Given an elementwise nonnegative vector $v \in \R^n_{>0}$, the $v$-weighted $L_1$ norm $\abs{\cdot}_{v}$ {of a vector $x \in \R^n$ } is $\abs{x}_{v} := {v}^\top \abs{x}$.  The dual norm $\abs{\cdot}_v^*$ to $\abs{\cdot}_v$ is the function $\abs{\cdot}_v^* := \max(|x|./v)$. In particular, the $L_1$ norm has an expression of $\norm{\cdot}_1 = \1^\top \abs{\cdot} = \abs{\cdot}_{\1}.$ {In particular, the dual norm to the $L_1$ norm is the $L_\infty$ norm $\norm{\cdot}_\infty = \norm{\cdot}_{1}^*$. The induced $L_\infty$ norm of a matrix $M \in \R^{n \times m}$ is $\norm{M}_{\infty} = \max_{i\in 1..n} \sum_{j=1}^m \abs{M_{ij}}$. 
}

{The $\ell_\infty$ and $\ell_1$ norms of a  sequence $y = \{y_t\}_{t \in \N}$ are
\begin{align}
    \norm{y}_\infty &= \text{sup}_{t\in\N} \norm{y_t}_\infty, & \norm{y}_1 &= \textstyle \sum_{t=0}^\infty \norm{y_t}_1 & 
\end{align}}

\subsection{Positive Linear Systems}
 
The {linear} discrete-time system under consideration has time axis $t \in \N$, state $x_t \in \R^n$, external input $w_t \in \R^e$, and measured output $y_t \in \R^p$:

\begin{align}
\label{eq:sys_xi_aut}
     x_{t+1} &= A x_t +  {E} w_t, \\ 
    y_t &= C x_t  + {F} w_t. \nonumber
\end{align}

The linear system in \eqref{eq:sys_xi_aut} is \textit{internally positive} if $x_0 \in \R^n_{\geq 0}$ and $\forall t \in \N: w_{t} \in \R^{{e}}_{\geq 0}$ implies $\forall t: y_{t} \in \R^p_{\geq 0}, x_t \in \R^n_{\geq 0}$ \cite{farina2000positive}. 
A necessary and sufficient condition for the state-space system in \eqref{eq:sys_xi_aut} to be \textit{internally positive} is if all matrices $(A,{E, C}, { F} )$ have nonnegative entries.  {Refer to \cite{farina2000positive} for a full review of positive systems theory}.

{The discrete-time system in \eqref{eq:sys_xi_aut} is Schur-stable if $A$ is a Schur matrix, implying that $\lim_{t\rightarrow \infty} x_t = 0$ for all $x_0 \in \R^n$ if $w_t = 0$ at all times $t \in \N$.} The following equivalent conditions ensure {that a matrix ${A} \in \R^{n \times n}_{\geq 0}$ is Schur} \cite{rantzer2018tutorial}:
\begin{enumerate}
\item The spectral radius of $A$ is less than 1.
\item There exists  $v_\infty \in \R^n_{>0}$ such that $v_\infty - A v_\infty \in \R^n_{>0}$.
\item There exists  $v_1 \in \R^n_{>0}$ such that $v_1 - A^\top v_1 \in \R^n_{>0}$.
\end{enumerate}

Conditions 2 and 3 refer to the existence of the following Lyapunov functions described by vectors $v_\infty, v_1 \in \R_{>0}^n$:
\begin{subequations}
\begin{align}
V_{{\infty}}(x) &= \max(x ./ v_\infty)  & 
V_{{1}}(x) &= v_1^\top x.
\end{align}
\end{subequations}

\section{Incremental Gains of \lure Systems}

{The interconnection of a linear system and a time-varying static nonlinearity $f(t, \cdot): \R^q \rightarrow \R^d$ is described as}
\begin{subequations}
\label{eq:sys_xi_lure}
\begin{align}
     x_{t+1} &= A x_t + B_1 z_t +   B_2 w_t \\
     z_t &= f(t, C_1 x_t + F_1 w_t) \\
    y_t &= C_2 x_t  + F_2 w_t.
\end{align}
\end{subequations}

The {\lure} system in  {\eqref{eq:sys_xi_lure}} has an $\ell_\infty$ incremental gain of $\eta$  if the following relationship holds for any input-output pair $(w^1, y^1), (w^2, y^2)$ of {\eqref{eq:sys_xi_aut} starting from the same initial condition $x_0$}:
\begin{align}
 \eta &= \sup_{w^1, w^2}  \frac{\norm{y^1 - y^2}_\infty}{\norm{w^1 - w^2}_\infty}. \label{eq:incr_gain}
    \intertext{The system in \eqref{eq:sys_xi_aut}  has an incremental $\ell_1$ gain of $\gamma$  if }
    \gamma &= \sup_{w^1, w^2} \lim_{T \rightarrow \infty} \frac{\sum_{t=0}^T \norm{y_t^1 - y_t^2}_1}{\sum_{t=0}^T \norm{w_t^1 - w_t^2}_1}. \label{eq:incr_gain_2}
\end{align}

This subsection will focus on deriving incremental gain bounds for systems in \eqref{eq:sys_xi_lure} via linear programs.

\label{sec:gain_compute}

\subsection{Assumptions}
Assumptions of system positivity and boundedness of the nonlinearity $f(t, \cdot)$ will now be detailed:

\begin{enumerate}    
    \item[A1] There exists a known matrix $\Delta \in \R^{{d \times q}}_{\geq 0}$ such that for all $ \zeta_1, \zeta_2 \in \R^q$,
    \begin{align}
         \sup_{t \in \N} \abs{f(t, \zeta_1) - f(t, \zeta_2)} &\leq \Delta \abs{\zeta_1 - \zeta_2}. \label{eq:incr_lure}
    \end{align}
    \label{assum:incr_lure}
    \item[A2] The matrices $(B_1, B_2, F_1, F_2)$ are nonnegative.
    \item[A3] The matrices $(A, C_1, C_2)$ are nonnegative.

\end{enumerate}

Assumption A1's statement in \eqref{eq:incr_lure} can be interpreted as an elementwise Lipschitz bound \cite{piengeon2024linear}, and allows for a linear-programming-based verification of stability and performance.
Assumption A2 and A3 ensure that the linear system in \eqref{eq:sys_xi_aut} is positive. {Assumption A3 will be removed in state-feedback controller synthesis.}

\subsection{Incremental  gain computation}

{Based on the \lure bound $\Delta$ in  Assumption A{1}, we define the following matrices:
\begin{align}
    A_\Delta &:= A + B_1 \Delta C_1, & B_\Delta&:= B_2 + B_1 \Delta F_1.
\end{align}}
We first note a linear programming method to bound the incremental $\ell_1$ gain from \cite{piengeon2024linear}, as modified for discrete-time ($h -A_\Delta{h}$) rather than continuous time  ($-A_\Delta{h}$):
\begin{prop}[Thm 4.1, Cor. 2.1 \cite{piengeon2024linear}]
\label{prop:l1_lure_analysis}
 The system in \eqref{eq:sys_xi_lure}
    has an $\ell_1$ incremental gain $w \rightarrow y$ upper-bounded by $\gamma$ {under} Assumptions  A1{-A3} if the following program is feasible:
    \begin{subequations}
\label{eq:l1_d_delta_analysis}
\begin{align}
\find_{h \in \R^n_{>0}} \quad &h -A_\Delta^\top h - C_2^\top\1_n  {\in \R^{n}_{>0}} \label{eq:l1_d_lure_dyn_analysis}\\
& \gamma \1_p - F_2^\top \1_e - B_\Delta^\top h {\in \R^{e}_{>0}}.\label{eq:l1_d_lure_gx_analysis}
\end{align}
\end{subequations}
\end{prop}

Our first contribution is a program to upper-bound the 
incremental $\ell_\infty$ gain of positive \lure systems in \eqref{eq:sys_xi_aut}.
\begin{thm}
\label{thm:linf_lure_analysis}
The affine system in  
\eqref{eq:sys_xi_lure} has an  incremental $\ell_\infty$  gain $w \rightarrow y$ upper-bounded by $\eta$ {under} Assumptions  A1{-A3} {if} the following constraints are feasible:
\begin{subequations}
\label{eq:linf_d_lure_analysis}
\begin{align}
\find_{v \in \R^n_{>0}} \quad &v -A_\Delta v - B_\Delta \1_e {\in \R^{n}_{>0}}\label{eq:linf_d_lure_dyn_analysis}\\
& \eta \1_p - F_2 \1_e - C_2 v {\in \R^{e}_{>0}}. \label{eq:linf_d_lure_gx_analysis}
\end{align}
\end{subequations}

\end{thm}
\begin{proof}
See Appendix \ref{app:linf_proof}.

\end{proof}

{We specifically highlight a state relation. Let $\abs{\delta x_{t}} = \abs{x^1_t - x^2_t}$ denote the time-wise absolute value of the increment between any two states in a trajectory of \eqref{eq:sys_xi_lure}, with similar incremental definitions for $\abs{\delta y_t}, \abs{\delta w_t}$. Let $\abs{\delta m_t} := \max_{t' \in 0..t} \abs{\delta w_t}$ denote the running maximum of $\abs{\delta w_t}.$ Given that strict feasibility of $v \in \R_{>0}^n$ implies that $v$ is non-infinite (bounded), and \eqref{eq:linf_d_lure_dyn_analysis} implies the existence of a  $\beta>0$ with  $v -A_\Delta v - B_\Delta \1_e -  \beta v {\in \R^{n}_{>0}}$, any trajectory of \eqref{eq:sys_xi_aut} satisfies the state relation
\begin{align*}
\abs{\delta x_{t+1}}_v^*- \abs{\delta x_t}_v^* \nonumber&\leq  \abs{B_\Delta \abs{\delta w_t}}_{v}^* - \norm{B_\Delta \abs{\delta x_t}}_\infty (1- \beta),
\end{align*}
and the output relation
\begin{align*}
    \norm{\delta y_t}_\infty &\leq \norm{C_2 \text{diag}(v)}_\infty \norm{\abs{\delta x_t ./ v} - \abs{\delta m_t}}_\infty + \eta \norm{\delta m_t}_\infty. \nonumber
    % \label{eq:gain_end_norm}
\end{align*}

The output $\norm{\delta y_t}_\infty $ will therefore converge to $\eta \norm{\delta m_t}_\infty$ at rate $(1-\beta)$. This decay to a robust invariant set has also been observed in the literature of superstable control systems \cite{polyak2002superstable} (linear systems where {$\norm{\cdot}_{\infty}$} is a polyhedral Lyapunov function).

The mathematical expression in \eqref{eq:linf_d_lure_analysis} is identical to formulations to compute the $\ell_\infty$-gain of a positive linear system with a matrix representation $(A + B_1 \Delta C_1, B_2 + B_1 \Delta F_1, C_2, F_2)$ (e.g. \cite[Lemma 16]{naghnaeian2014linfinity}), which would classically require these four matrices to have all nonnegative entries. The assumption in A1 that $(A, B_1, C_1)$ are all positive is required to apply the $\Delta$-bound from assumption A2 to when describing  \lure uncertainty in \eqref{eq:incr_lure}.}

{The incremental $\ell_{\infty}$ norm is the dual of  the incremental $\ell_{1}$ norm. The program to compute the incremental $\ell_\infty$ norm with matrices $(A + B_1 \Delta C_1, B_2 + B_1 \Delta F_1, C_2, F_2)$  in \eqref{eq:linf_d_lure_analysis} is the same as computing the incremental $\ell_1$  of the dual system $(A_\Delta^\top, C_2^\top, B_\Delta^\top,  F_2^\top)$. However, this duality of linear systems can break down when considering the \lure connection in \eqref{eq:sys_xi_lure}. Instead, an adjoint system must be proven to be backward-in-time stable, as detailed in Proposition 2.13 of \cite{briat2021hybrid}. Our direct proof in Appendix \ref{app:linf_proof} applies the elementwise-Lipschitz bound $\Delta$ in \eqref{eq:incr_gain} without requiring adjoint system theory.}

\section{Incremental Gain-Bounded Synthesis}
\label{sec:regulate}

We develop full-state-feedback controllers in the setting of elementwise uncertainty  such that the closed-loop system has bounded incremental {$\ell_\infty$} performance.

The system model in \eqref{eq:sys_xi_lure} is augmented with a controlled input $u \in \R^m$ to form the system
\begin{subequations}
\label{eq:sys_xi_lure_u}
\begin{align}
     x_{t+1} &= A x_t + B_1 z_t +   B_2 w_t  + B_3 u_t\\
     z_t &= f(t, C_1 x_t + F_1 w_t + D_1 u_t) \\
    y_t &= C_2 x_t  + F_2 w_t + D_2 u_t.
\end{align}
\end{subequations}

The class of state-feedback policies considered in this section are affine controllers parametized by a matrix $K \in \R^{m \times n}$ and a vector $g \in \R^m$ as {$u_t = K x_t + g$.}
The affine offset $g$ reflects a sustained input to the system, and is assumed to be given (or optimized in a previous layer of a hierarchical control system). The affine offset is irrelevant in the incremental gain analysis, because the same offset is applied to any two trajectories of the controlled system.

The state-feedback control problem to regulate the incremental $\ell_\infty$ performance:
\begin{prob}
    Given a system \eqref{eq:sys_xi_lure_u} and an affine offset $g$, find a matrix $K$ such that the law ${u=Kx + g}$ minimizes the {incremental} $\ell_\infty$ gain {from} $w$ and $y$.
\end{prob}

\begin{thm}
\label{thm:linf_regulate}
    A sufficient condition for the existence of  a matrix $K$ such that the law ${u=Kx + g}$ regulates \eqref{eq:sys_xi_lure_u} with an incremental $\ell_\infty$ gain at most $\eta$ is if the following program is feasible (under Assumptions A{1}-A3 and given the affine offset $g$):
        % \begin{subequations}
\begin{align}
\find_{v, X,   Y} \quad &v -(A_{\Delta}X  + {(B_3 + B_1 \Delta D_1)}) Y) \1_n - B_\Delta \1_e  {\in \R^{e}_{>0}} \\
& \eta \1_p - F_2 \1_e - C_2v - D_2 Y \1_n {\in \R^{n}_{>0}}\label{eq:linf_d_lure_control}\\
& A X + B_3 Y \in \R^{n \times n}_{\geq 0}, \nonumber \
 C_1 X + D_1 Y \in \R^{q \times n}_{\geq 0}, \nonumber \\
& C_2 X + D_2 Y \in \R^{p \times n}_{\geq 0}, \nonumber \\
& v \in \R^n_{>0}, \quad X = \diag{v}, \quad  Y \in \R^{m \times n}. \nonumber
\end{align}
% \end{subequations}
If program \eqref{eq:linf_d_lure_control} is feasible, then the controller $K$ can be recovered by computing $K = Y X^{-1}$.
\end{thm}
\begin{proof}
See Appendix \ref{app:linf_proof_control}
\end{proof}

\section{Numerical Examples}
\label{sec:examples}

MATLAB 2023b code to generate all experiments is publicly available \footnote{\url{https://github.com/Jarmill/pos_lure_lp}}. Dependencies include YALMIP \cite{lofberg2004yalmip} and Mosek \cite{mosek92}.

\subsection{Leslie Model Analysis}

The first example involves a Leslie model for population dynamics \cite{leslie1945use} {\cite{cushing1998introduction}}. This system has $n=5$ states, $e=2$ external inputs, $q=2$ \lure outputs, and $p=1$ output channel. The matrices that describe this system are:
\begin{align}
    A &= \begin{bmatrix}
        0 & 0.1 & 0.6 & 0.2 & 0.1 \\
        0.95 & 0 & 0 & 0 & 0 \\
        0 & 0.9 & 0  & 0 & 0 \\
        0 & 0 & 0.7 & 0 & 0 \\
        0 & 0 & 0 & 0.5 & 0 \\
    \end{bmatrix}, \nonumber
    \end{align}
    \begin{align}
    B_1 &= \begin{bmatrix}
        I_2 \\ \0_{3 \times 2}
    \end{bmatrix}, & B_2 &= \begin{bmatrix}
        I_2 \\ \0_{3 \times 2}
    \end{bmatrix}, & F_1 &= \begin{bmatrix}
        0.1 & 0.1 \\ 0.1 & 0.1
    \end{bmatrix} \nonumber\\
    C_1 &= \begin{bmatrix}
        0 & 0 & 1 & 1 & 0 \\ 0 & 0 & 0 & 1 & 1
    \end{bmatrix}, & 
    C_2 &= \begin{bmatrix}
        0 & 0 & 0 & 0 & 1 
    \end{bmatrix},  & F_2 &= \begin{bmatrix}
        0 & 0
    \end{bmatrix}. \label{eq:leslie_model}
\end{align}

\subsubsection{Element{wise}-Bounded Nonlinearity}

We first consider the circumstance where the \lure nonlinearity is bounded according to the elementwise Assumption A2 with parameter $\Delta = {\tau} I$ for a scalar value ${\tau} \geq 0$. {The} specific \lure nonlinearity used in this example is 
\begin{align}
    f(t, \zeta) &= {\begin{bmatrix}
        0.5 \\ 0
    \end{bmatrix} + \frac{\tau}{2} \begin{bmatrix}
        \zeta_1 + \sin(\zeta_1)\\ \zeta_2 + \sin(\zeta_2)
    \end{bmatrix}}. \label{eq:leslie_lure}    
\end{align}

In the absence of an applied input ($w=0$), the fixed point $x^*$ of dynamics \eqref{eq:sys_xi_lure} with the plant matrices \eqref{eq:leslie_model} and the nonlinearity \eqref{eq:leslie_model} is
\begin{equation}
    x^* = \begin{bmatrix}
               2.6505 &
    2.5952&
    2.3356&
    1.6349&
    0.8175
    \end{bmatrix}^\top.
\end{equation}

\begin{figure}[h]
    \centering
    \includegraphics[width=0.5\linewidth]{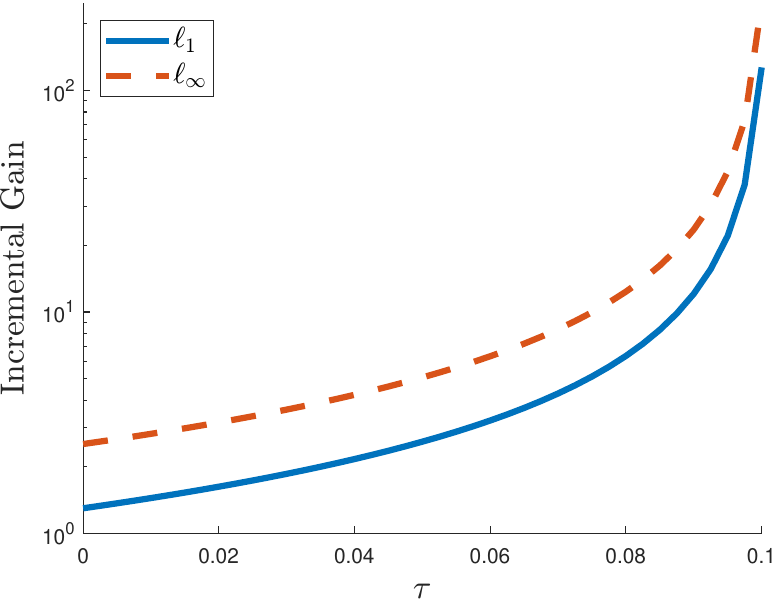}
    \caption{Upper-bounds on the incremental $\ell_1$ and $\ell_\infty$ gains {of the Leslie system in \eqref{eq:leslie_model} as the bound $\tau$ increases}}
    \label{fig:leslie_l1_linf}
\end{figure}

Figure \ref{fig:leslie_l1_linf} reports upper bounds on the incremental gains $\ell_1$ and $\ell_\infty$ {from} $w \rightarrow y$ for the Leslie system. All programs become infeasible at $\tau \geq 0.125$. 
{We are therefore unable to certify absolute stability  of the \lure system nor finiteness of the $\ell_1$ or $\ell_\infty$ incremental gains when $\tau \geq 0.125.$ 
}
% and the reported upper bounds on the {incremental} gain{s are} thus $\infty$.

\begin{figure}[!h]
    \centering
    \includegraphics[width=0.5\linewidth]{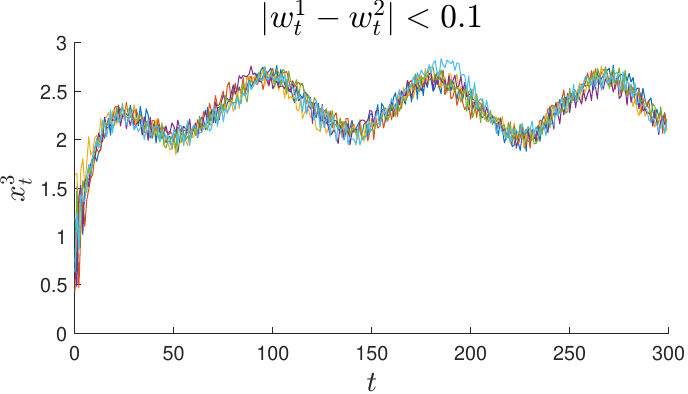}
    \caption{Simulation of the state $x^3_t$ in the interconnection of \eqref{eq:leslie_model} and  \eqref{eq:leslie_lure} {starting at $x_0 = \1_5$ as perturbed by an $\ell_\infty$-bounded signal $w$ drawn from  \eqref{eq:w_sample}}}
    \label{fig:leslie_x}
\end{figure}

{

{We perturb the Leslie \lure} system with an applied chirp signal plus $\ell_\infty$-bounded randomly sampled noise:
\begin{align}
    w(t) = 0.15  \cos\left(\frac{\pi}{50}t  +  10^{-4} t^2\right) + w^{\text{rand}}(t), \label{eq:w_sample}
\end{align}
in which $w^{\text{rand}}(t)$  is uniformly drawn from a uniform distribution between [-0.05, 0.05] at each discrete time $t$.
Figure \ref{fig:leslie_x} plots a trace of the  third coordinate $x^3(t)$ {in the} Leslie system under the  \lure input in \eqref{eq:leslie_lure} 
with parameter ${\tau} = 0.05$ and an input signal in {\eqref{eq:w_sample}}.
 The maximum $\ell_\infty$ distance between any two $w$ samples in \eqref{eq:w_sample} is $0.1$. The corresponding $\ell_\infty$ gain is 5.0665.
}

\begin{figure}[h]
    \centering
    \includegraphics[width=0.5\linewidth]{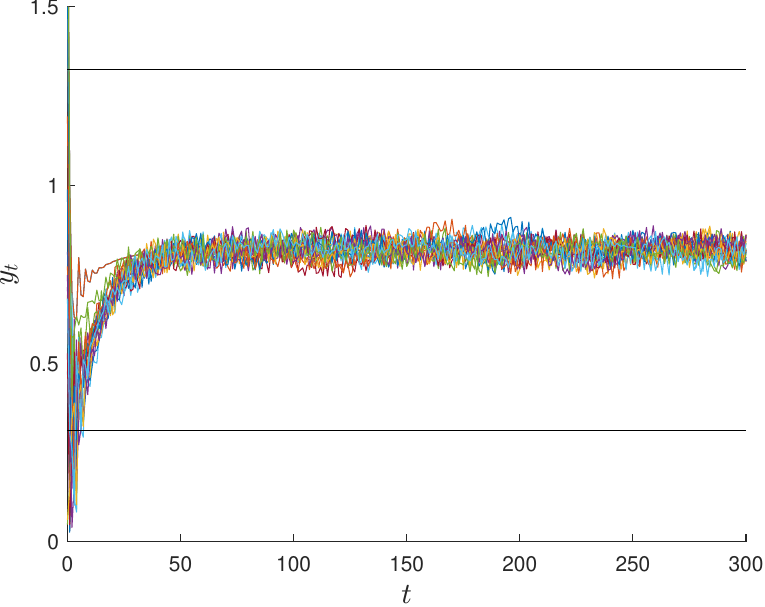}
    \caption{{Sampled outputs  $y$ and a certified incremental $\ell_\infty$ bound for  $\norm{y_t - y^*}_\infty$ of the interconnection of \eqref{eq:leslie_model} and \eqref{eq:leslie_lure} w.r.t. $\delta = 0.05$ and input $\norm{w^1-w^2}_\infty < 0.1$ from  \eqref{eq:w_sample}.}}
    \label{fig:leslie_y}
\end{figure}

Figure \ref{fig:leslie_y} displays the bounds of $0.1*5.0665 = {0.50665}$ in black lines above and below the unperturbed fixed point of $y^* = 0.8175$, {along with} 20 simulated trajectories. After transients, the trajectories stay within the region bordered by the black lines.

\subsection{System Control}
\label{sec:system_control}
The second example involves a randomly generated positive linear system with $n=10$ states, $m=4$ controlled inputs, $e=2$ external inputs, $d=2$ \lure inputs, and  $p=2$ outputs. 

\begin{figure}[h]
    \centering
    \includegraphics[width=0.5\linewidth]{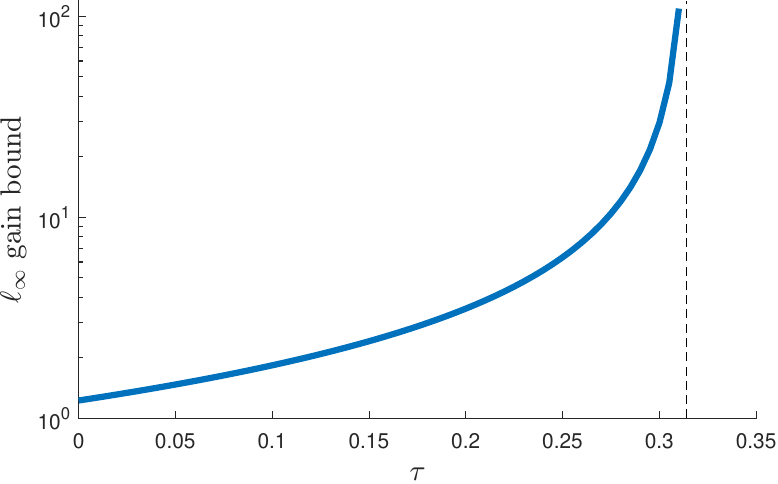}
    \caption{{Achieved upper-bounds on the incremental $\ell_1$ and  $\ell_\infty$  gains by state-feedback-regulation for the randomly generated system in Section \ref{sec:system_control} as ${\tau}$ increases.}}
    \label{fig:regulate}
\end{figure}

The uncontrolled system 
% in \eqref{eq:sys_second} 
is unstable with a spectral radius of $\rho(A) = 1.1888$. Figure \ref{fig:regulate} {plots upper-bound on the incremental $\ell_\infty$ norm of the closed-loop system,} as obtained by solving the Linear Program in \eqref{eq:linf_d_lure_control} versus a \lure nonlinearity described by the incremental bound $\Delta = {\tau} I_2$ as ${\tau}$ increases. The maximal possible ${\tau}$ value that can be certifiably regulated via Program \eqref{eq:linf_d_lure_control} (up to five decimal places) is $0.31373$.

\section{Conclusion}
\label{sec:conclusion}

This work presented a linear programming formulation for upper-bounding the incremental  $\ell_\infty$ gains of a positive linear system under \lure uncertainty. 
Linear programming can also be used to regulate incremental $\ell_\infty$ performance using state-feedback synthesis in a tractable manner.
Future work includes {development of output feedback controllers for incremental stabilization,} and  broadening the set of uncertainties that can be handled through copositive techniques (e.g. incremental positive $\ell_2$ gains \cite{ebiharal2+}). A challenging but promising area for future research is in treating automated synthesis of optimization algorithms \cite{scherer2023optimization} under positivity constraints.

\section*{Acknowledgements}

% The author would like to thank Prof. Carsten Scherer, Sophie Hall, Sarah H. Q. Li, Jaap Eising, and Zhiyu He for discussions about positive systems and incremental gains.
This work was funded by Deutsche Forschungsgemeinschaft (DFG, German Research Foundation) under Germany's Excellence Strategy - EXC 2075 – 390740016. We acknowledge the support by the Stuttgart Center for Simulation Science (SimTech).

\appendix

\section{Proof of Incremental Gains}
\label{app:proofs}
This appendix contains the proofs of Theorems \ref{thm:linf_lure_analysis} and \ref{thm:linf_regulate}.
% \subsection{Technical Lemmas}
We first note basic properties of nonnegative matrices. If  $c, z \in \R^n$ are vectors, then $\abs{c^\top z} \leq \abs{c}^\top \abs{z}$. 
% \begin{lem}
   Also, for {all} matrices $M \in \R^{p_1 \times p_2}_{\geq 0}$, $\Delta \in \R^{p_2 \times p_3}_{\geq 0}$, $N \in \R^{p_3 \times p_4}_{\geq 0}$, and $\Sigma \in \R^{p_2 \times p_3}_{\geq 0}$; {it holds that  $\Sigma \leq \Delta$ implies $M \Sigma N \leq M \Delta N$.}

\subsection{Proof of Theorem \ref{thm:linf_lure_analysis}}
\label{app:linf_proof}
% \begin{proof}

{Let} $(x^1, x^2) \in (\R^n_{\geq 0})^2$ be a pair of states, $(w^1, w^2) \in (\R^e)^2$ be a pair of disturbances, and $(z^1, z^2) \in (\R^d)^2$ be a pair of \lure outputs from $f$. The incremental variable for the state $x$ is $\delta x = x^1 - x^2$ (with similar incremental definitions {performed} for the input $\delta w$, measured output $\delta y$, and \lure output $\delta z$).

{
We consider an incremental storage function $V(\delta x_t) = \norm{\delta x_t ./ v}_{\infty} = \abs{x_t}_v^*.$ The change of this storage function is bounded by
\begin{subequations}
\begin{align*}
\delta V_t &:= {V(\delta x_{t+1}) - V(\delta x_t)} \\
     &= {\abs{\delta x_{t+1}}_v^* - \abs{\delta x_t}_v^*} \\
    &= \abs{A \delta x_t + B_1 \delta z_t + B_2 \delta w_t}_{v}^*- \abs{\delta x_t}^*_{v}. \\
    \intertext{We apply absolute values to all incremental terms to produce an upper-bound}
    &\leq  \abs{(A \abs{\delta x_t} + B_1 \abs{\delta z_t} + B_2 \abs{\delta w_t})}_{v}^* - \abs{\delta x_t}_v^*. \\
     \intertext{Defining $A_\Delta:=A+B_1 \Delta C$ and $B_\Delta:=B_2 + B_1 \Delta F_1$ using the matrix $\Delta$ from Assumption A1    {and using the equality $\abs{(\abs{\delta x_t})}^*_v = \abs{\delta x_t}^*_v,$}
     the $\abs{\delta z_t}$ term can be bounded as}
     &\leq  \abs{(A_\Delta \abs{\delta x_t} + B_\Delta \abs{\delta w_t})}_{v}^* - \abs{\delta x_t}_v^*. \\
     &\leq  (\abs{A_\Delta \abs{\delta x_t}}_{v}^*  - \abs{\delta x_t}_v^*) + \abs{B_\Delta \abs{\delta w_t}}_{v}^*
\end{align*}
\end{subequations}
A consequence of \eqref{eq:linf_d_lure_dyn_analysis} is that there exists a 
$\beta > 0$ with 
${v -A_\Delta v - B_\Delta \1_e -\beta \1_n \in \R^n_{>0},}$
{which implies that}
\begin{align*}
   ( \abs{\delta x_t}_v^* - \abs{A_\Delta \abs{\delta x_t}}_{v}^*) \geq  \norm{B_\Delta \abs{\delta x_t}}_\infty    + \beta \abs{\delta x_t}_v^*
\end{align*}
Applying this relation to the storage function decrease term yields
   \begin{align*}
   \delta V_t &\leq  \abs{B_\Delta \abs{\delta w_t}}_{v}^* - \norm{B_\Delta \abs{\delta x_t}}_\infty - \beta \abs{\delta x_t}_v^* 
   % \label{eq:beta_decay}
\end{align*}
We now bound the output coordinates of the incremental output as
\begin{align*}
    {\delta y_t} &= {C_2 \delta x_t + F_2 \delta w_t} \\
    \abs{\delta y_t} &\leq C_2 \abs{\delta x_t} + F_2 \abs{\delta w_t} \\
    &= (C_2 \text{diag}(v)) \abs{\delta x_t ./ v} + F_2 \abs{\delta w_t} 
    \intertext{Defining the running incremental maximum $\abs{\delta m_t} := \max_{s=0..t} \abs{\delta w_s}$, the incremental output bound satisfies }
   \abs{\delta y_t} &\leq (C_2 \text{diag}(v)) \abs{\delta x_t ./ v} + F_2 \abs{\delta m_t}.
   \intertext{Feasibility of \eqref{eq:linf_d_lure_gx_analysis} allows for a substitution of the $F_2 \abs{\delta m_t}$ term, leading to   }
   \abs{\delta y_t} &\leq (C_2 \text{diag}(v)) (\abs{\delta x_t ./ v} - \abs{\delta m_t}) + \eta \abs{\delta m_t}    \\
      \norm{\delta y_t}_\infty &\leq \norm{(C_2 \text{diag}(v)) (\abs{\delta x_t ./ v} - \abs{\delta m_t})}_\infty + \eta \norm{\delta m_t}_\infty 
\end{align*}
}
The $\eta$ factor between  $\norm{\delta m_t}_\infty$ and $\norm{\delta y_t}_\infty$ proves the desired $\ell_\infty$ incremental gain relation.$\qed$
% \end{proof}
\subsection{Proof of Theorem \ref{thm:linf_regulate}}
\label{app:linf_proof_control}
% \begin{proof}
    Substitution of the state-feedback policy $u = K x + g$ into \eqref{eq:sys_xi_lure_u} leads to the $(u, x)$ closed-loop expression for dynamics
    % \begin{subequations}
\begin{align}
     x_{t+1} &= (A + B_3 K) x_t + B_1 z_t +   B_2 w_t  + B_3 g \nonumber \\
     z_t &= f(t, (C_1 + D_1 K) x_t + F_1 w_t + D_1 g) \label{eq:sys_xi_lure_u_closed}\\
    y_t &= (C_2 + D_2 K) x_t  + F_2 w_t + D_2 g. \nonumber
\end{align}
The terms $(B_3 g, D_1 g, D_2g)$ are all constants among all trajectories of \eqref{eq:sys_xi_lure_u_closed}, and can therefore be discarded from an incremental analysis. Application of Theorem \ref{thm:linf_lure_analysis} to the system in \eqref{eq:sys_xi_lure_u_closed} ensures that satisfaction of the following constraints will ensure that the law $u=Kx + g$ respects an incremental gain $\ell_\infty$ less than $\eta$:
        \begin{subequations}
\label{eq:linf_d_lure_control_K}
\begin{align}
\find_{v, K} \quad &v -(A_\Delta + (B_3 + B_1 \Delta D_1) K) v\nonumber  - B_\Delta \1_e {>0} \\
& \eta \1_p - F_2 \1_e - (C_2 + D_2 K) v {>0}\label{eq:linf_d_lure_gx_control_K}\\
 & C_1 + D_1 K \in \R^{n \times e}_{\geq 0}, \ \label{eq:linf_d_lure_gx_control_K_pos1} C_2 + D_2 K \in \R^{p \times n}_{\geq 0}  \\
& A  + B_3 K \in \R^{n \times n}_{\geq 0},\label{eq:linf_d_lure_gx_control_K_pos3} \\
&K \in \R^{m \times n}_{\geq 0}, v \in \R^n_{>0}. \nonumber
\end{align}
\end{subequations}
Constraints \eqref{eq:linf_d_lure_gx_control_K_pos1}-\eqref{eq:linf_d_lure_gx_control_K_pos3} ensure that Assumption {A3} is respected for the closed-loop system \eqref{eq:sys_xi_lure_u_closed}.
We {now}  deploy a standard convexification \cite{boyd1994linear} by defining a new variable $Y$ 
such that $K \diag{v} = Y$ to break the $K$-$v$ multiplication. Such a $Y$-redefinition was also performed in the positive system setting for $\ell_\infty$ regulation in Lemma 2 of \cite{briat2013positive} (continuous-time) and  Lemma 16 of \cite{naghnaeian2017performance} (discrete-time). 

Right-multiplication by $(A+B_3 K)$ by the positive-diagonal matrix $\diag{v}$ does not change the its nonnegativity properties, thus proving the theorem.

\bibliographystyle{ieeetr}        % Include this if you use bibtex
\bibliography{hybrid}

\end{document}